\documentclass[11pt, reqno]{amsart}
\usepackage{amsmath, amsthm, amscd, amsfonts, amssymb, graphicx, color}
\usepackage[bookmarksnumbered, colorlinks, plainpages]{hyperref}
\usepackage{enumitem}

\makeatletter \oddsidemargin.9375in \evensidemargin \oddsidemargin
\def\authorsaddresses#1{\dedicatory{#1}}
\newtheorem{theorem}{Theorem}[section]
\newtheorem{lemma}[theorem]{Lemma}
\newtheorem{proposition}[theorem]{Proposition}

\theoremstyle{definition}
\newtheorem{definition}[theorem]{Definition}

\theoremstyle{remark}

\theoremstyle{approach}

\numberwithin{equation}{section}

\begin{document}
\setcounter{page}{1}


\title[TOPOLOGICAL INVARIANT MEANS ON $B^*_p(G)$ AND WEAKLY COMPACT Maultipliers of $A_p(G)$ ]{TOPOLOGICAL INVARIANT MEANS ON $B^*_p(G)$ AND WEAKLY COMPACT Maultipliers of $A_p(G)$}

\author[H. G. AMINI AND A.REJALI]{H. G. AMINI AND A.REJALI}


\begin{abstract}
In this paper we investigate the compact and weakly compact multipliers of the Herz-algebras $A_p(G)$. Let $B_p(G)$ be the space of pointwise multipliers of $A_p(G)$. We show that there is a topological invariant mean on $B^*_p (G)$. Furthermore, we show that if $B^*_p(G)$ is separable, then $G$ is a discrete group.
\end{abstract}

\maketitle


\section{Introduction}
Let $G$ be a locally compact group. For $1 < p < \infty$, let $A_p(G)$ denote the Herz algebra of $G$ and $B_p(G)$ denote the space of multipliers of $A_p(G)$. The results obtained in this paper improve some results obtained by E. E. Granirer \cite{6} for $A^*_p(G)$. Also, we give a simple new proof for Proposition 7 of \cite{6}.\\
\indent
In section 3 we both define and show the existence of topological invariant means on $B^*_p(G)$. Also we show that the invariant mean is necessarily unique on $wap(B_p(G))$, the space of weakly almost periodic functional on $B_p(G)$. Furthermore, if $B^*_p(G)$ is separable then $G$ is discrete. \\
\indent
In section 4 we show that if $A_p(G)$ has a [weakly] compact multiplier, then $G$ is discrete. Furthermore, if $G$ is amenable, then each [weakly] compact multiplier is in $A_p(G)$.
\section{Preliminaries and some notations}
Let $G$ be a locally compact group equipped with a fixed left Haar measure $\lambda= dx$. The spaces $L^p(G)$ ($1\leq p\leq\infty$) have their usual meaning. Let $C_{00}(G)$, denote the space of complex continuous functions on $G$ with compact support. Also, let $C_0(G)$ be the space of complex continuous functions on $G$ which tend to 0 at infinity with $\|.\|_{\infty}$-norm. Then $M(G) = C_0(G)^*$  denotes the bounded Borel complex measures	on $G$ with convolution as multiplication, and with variation norm. For $f$ in $L^1(G)$,
the operator $\rho(f): L^p(G) \longrightarrow L^p(G)$ is defined by $\rho(f)(g) = f*g$. Clearly $\rho(f)$ is a bounded linear operator on $L^p(G)$ and by $\|\rho(f)\|$ we shall always denote the operator norm of $\rho(f)$. \\
\indent
For $1 < p < \infty$ let $A_p(G)$ denote the linear subspace of $C_0(G)$ consisting of all
functions of the form $u(x)=\sum_{i=1}^{\infty}g_i*f^{\vee}_i(x)$ where $f_i\in L_p(G)$, $g_i\in L_q(G)$, $\frac{1}{p}+\frac{1}{q}=1$, $\sum_{i=1}^{\infty}\|f_i\|_p\|g_i\|_q<\infty$, $f^{\vee}(x)=f(x^{-1})$ for $x\in G$. $A_p(G)$ is a commutative Banach algebra with respect to pointwise multiplication and the norm,
$$\|u\|_{A_p(G)}=\inf\{\sum_{i=1}^{\infty}\|f_i\|_p\|g_i\|_q: u=\sum_{i=1}^{\infty}g_i*f^{\vee}_i\}.$$
In the case $p = 2$, $A_2(G) = A(G)$ is the Fourier algebra of $G$, which was introduced for non commutative groups by Eymard \cite{3}. For $p\not= 2$, $A_p(G)$ was first studied by Herz in \cite{7}.\\
\indent
We denote by $B_p(G)$ the set of bounded complex continuous functions $u$ on $G$ such that $uv\in A_p(G)$, for all $v\in A_p(G)$. The norm in $B_p(G)$ is given by
$$\|u\|=\sup\{\|uv\|_{A_p(G)}: \|v\|_{A_p(G)}\leq 1\}.$$
Then $B_p(G)$ is the space of pointwise multipliers of $A_p(G)$. Furthermore $B_p(G)$ is
a commutative Banach algebra with respect to pointwise multiplication. \\
\indent
Each element $f\in L^1(G)$ defines a bounded functional $\phi_f$ on $A_p(G)$ by
$$<\phi_f,u>=\int_Gf(x)u(x)dx,\qquad u\in A_p(G).$$
The norm of $\phi_f$ as an element of $A_p(G)^*$ and the operator norm of $\rho(f)$ are the
same. That is,
$$\|\phi_f\|=\sup_{\|u\|_{A_p(G)}\leq 1}|<\phi_f,u>|=\|\rho(f)\|=\sup_{\|g\|_p\leq 1}\|f*g\|.$$
It follows that $L^1(G)$ can be considered as a subspace of $A_p(G)^*$ \cite{1}. By definition, $PM_p(G)$ denote the closure of $L^1(G)$, considered as an algebra of convolution operators on $L^p(G)$, with respect to the weak operator topology in $B(L^p(G))$, the bounded operator on $L^p(G)$. The space $PM_p(G)$ can be identified with the dual of $A_p(G)$ for each $1 < p < \infty$ \cite{1}. If $G$ is Abelian with dual $\hat{G}$, then $A_2(G) = A(G) = L^1(\hat{G})$, $B_2(G) = B(G) = M(\hat{G})$ and $PM_p(G) = L^{\infty}(\hat{G})$. We define the module action of $B_p(G)$ on $B^*_p(G)$ by $<\phi u, v>=<\phi ,uv>$ for $\phi\in B^*_p(G)$, $u, v \in B_p(G)$. For $x\in G$ we show by $\chi_x=1_{\{x\}}$ the characteristic function of $\{x\}$.
\section{Invariant Means on $B^*_p(G)$}
Let $I:B_p(G)\longrightarrow C, I(u) = u(e)$. Denote
$$MB^*_p(G) =\{F\in B^*_p(G) : \| F\|= F(I) = 1\}$$
and for each $x\in G$,
$$S^p_A(x) =\{u \in A_p(G) : u(x) = \|u\|= 1\}$$
and
$$S^p_B(x) = \{u \in B_p(G) : u(x) = \|u\|= 1\}.$$
Then $S^p_A(x)$ and $S^p_B(x)$ are commutative semigroup and $S^p_A(x)\subseteq S^p_B(x)$. For $x = e$, denote $S^p_A = S^p_A(e)$ and $S^p_B = S^p_B(e)$. Also if $\pi: B_p(G)\longrightarrow B^{**}_p(G)$ is the canonical map it is clear that $\pi(S^p_B)\subseteq M B^*_p(G)$.\\
\indent
The proof of the next theorem is similar to the proof of \cite[Proposition 4]{6} and
by this theorem we can define the set of topological invariant means on $B^*_p(G)$.
\begin{theorem}
There is $F\in MB^*_p(G)$ such that, $uF = u(e)F$, for $u\in B_p(G)$.
\end{theorem}
\begin{definition}
$$TIMB^*_p(G) :=\{F\in MB^*_p(G) : uF = u(e)F,\, \text{for all $u\in B_p(G)$}\}$$ is called the
set of topological invariant means on $B^*_p(G)$.
\end{definition}
\begin{proposition}
$TIMB^*_p(G) \cap B_p(G)\not=\emptyset$ if and only if $G$ is discrete. In this case $TIMB^*_p(G) \cap B_p(G)=\{\chi_e\}$.
\end{proposition}
{\bf Proof.}
Let $G$ be discrete. Then it is clear that 
$\chi_e \in TIMB^*_p(G)\cap B_p(G)$. Conversely, if $u\in TIMB^*_p(G) \cap B_p(G)$. Then $u(e)=1$, and for each $v\in S^p_B$, $uv=u$. Let
$x\not= e$. Then there is $v\in S^p_B(G)$, $v(x) = 0$. So, $u(x) = uv(x) = 0$. Hence $u = \chi_e$ and $G$ is discrete. This show that $TIMB^*_p(G) \cap B_p(G)=\{\chi_e\}$. \\
 
The next theorem is the main result of this section. When $G$ is a discrete commutative group, Forrest show that if $PM_p(G)$ is separable, then $G$ is finite (see \cite{4} Theorem 1).
\begin{theorem}
Let $G$ be a locally compact group and $B^*_p(G)$ be separable.Then $G$ is discrete.
\end{theorem}
{\bf Proof.} Let $u$ be in $S^p_A$ with compact support. Let $E = supp u$ and $K = uS^p_A$. Then by the Markov- Kakutani fixed point theorem (see \cite{2}, p. 456) there is $F\in \{w^*clK\}\cap TIMB^*_p(G)$. Since $w^*$-topology on the unit ball of $B^{**}_p(G)$ is metrizable. Hence there exists a sequence $u_n\in S^p_A$ such that for all $\phi\in B^*_p(G)$, $\phi(uu_n)\longrightarrow F(\phi)$. \\
\indent
Now let $k\in A_p(G)$ with compact support such that $k = 1$ on $E$ and $\psi\in PM_p(G)$.
Then $k\psi$ is a bounded linear functional on $B_p(G)$ and $\psi(uu_n) =k\psi (uu_n)\longrightarrow F(k\psi)$. Therefore $(uu_n)$ is a weakly Cauchy sequence in
$$A^p_E(G) = \{v\in A_p(G) : supp v\subseteq E\}.$$
But by Lemma 18 on p. 131 of \cite{5} it follows that $A^p_E(G)$ is weakly sequentially complete for any compact set $E\subseteq G$. So there is $v_0\in A^p_E(G)$ such that for each $\psi\in PM_p(G)$, $\psi(uu_n)\longrightarrow \psi(v_0)$. Hence $\psi(v_0) = F(k\psi)$. Thus for each $w\in S^p_A(G)$ we have $\psi(v_0w) = F(k\psi w) = F(k\psi) =\psi(v_0)$. So, $v_0w = v_0$. Let $x\not= e$ and $w\in S^p_A(G)$ with $w(x) = 0$. We get $v_0(x) = 0$. However $v_0(e)=I(v_0) = F(kI) = F(I) = 1$. So $v_0 = \chi_e$. Since $v_0\in A_p(G)$ is continuous, so $G$ is discrete.
\begin{theorem}
Let $G$ be a locally compact group. Let $J$ be a closed ideal in $B_p(G)$ such that there is $u\in J$ with $u(e)\not= 0$. Then $TIMB^*_p(G)\subseteq J^{\perp\perp}$.
\end{theorem}
{\bf Proof}
Let $F\in TIMB^*_p(G)$ and $T\in J^{\perp}$. Since $uT = 0$. It follows that $u(e)F(T) = F(uT)= 0$. Hence $F(T) = 0$. Therefore $F\in J^{\perp\perp}$.\\

Let $\mu\in M(G)$ and 
$$< u,\mu>=\int_G u(x)d\mu(x)\qquad \text{for $u\in B_p(G)$}$$
then
$$|<u,\mu>|\leq \|\mu\| \|u\|_{\infty}\leq \|\mu\| \|u\|_{B_p(G)}.$$
Hence $\mu$ is a bounded linear functional on $B_p(G)$. Let $wap(B_p(G))$ denote the linear
space of all weakly almost periodic functional on $B_p(G)$. \\
\indent
Similar to the next proposition we can given a new proof for the Proposition 7 of \cite{6}.
\begin{proposition}
$M(G)\subseteq wap(B_p(G))$.
\end{proposition}
{\bf Proof.}
Since $C(G)$, the linear space of bounded continuous function on $G$ is a $C^*$-algebra and the unit ball of $B_p(G)$ is contained in the unit ball of $C(G)$. Now let $\mu\in M(G)$ and $(f_n)$ and $(g_m)$ be two sequence of the unit ball of $B_p(G)$ such that
$$\lim_n\lim_m <\mu, f_ng_m >, \qquad \lim_m\lim_n <\mu, f_ng_m>$$
both exist. Then by Arens regularity of $C(G)$, they are equal. So, $\mu$ is a weakly almost periodic functional on $B_p(G)$.\\
\indent
The proof of the next proposition is similar to the Proposition 9 \cite{6}.
\begin{proposition}
There is a unique topological invariant mean on $wap(B_p(G))$. In fact there is a unique $F\in wap(B_p(G))^*$ such that $F(I) = 1$ and $F(uT) = u(e)F(T)$, for all $u\in B_p(G)$ and $T\in wap(B_p(G))$.
\end{proposition}
\section{Compact and Weakly Compact Multipliers on $A_p(G)$}
In this section we shall present some results about compact and weakly compact multiplier of $A_p(G)$. For $u\in B_p(G)$, let $\Gamma_u : A_p(G)\longrightarrow A_p(G)$, defined by $\Gamma_u(v) = uv$ for $v\in A_p(G)$.
\begin{lemma}
If for some $u\not= 0$ in $B_p(G)$, $\Gamma_u$ is a weakly [compact] operator, then $G$ is discrete.
\end{lemma}
{\bf Proof.}
Let $u(x)\not= 0$ for some $x\in G$ and $K= w-cl\{uv : v\in S^p_A(x)\}$. Then $K$ is a weakly compact convex subset of $A_p(G)$ and For each $w\in S^p_A(x)$, $\Gamma_w(K)\subseteq K$. The operators $\Gamma_w$ are pairwise commuting and ($A_p(G)$, weak) to ($A_p(G)$, weak) continuous. By the Markov- Kakutani fixed point theorem (\cite{2}, p.456) there is $u_0\in K$ such that for each $w\in S^pA(x)$, $wu_0 = u_0$. Let $y\not= x$. Then there is $v\in S^p_A(x)$ with $v(y) = 0$. Hence $u_0(y) = v(y)u_0(y) = 0$. Therefore $u_0 = u(x)\chi_x$ and $G$ is discrete.\\

Now, we define the set of weakly compact multipliers on $A_p(G)$ as
$$WCB_p(G) := \{u\in B_p(G): \Gamma_u\,  \text{is weakly compact operator on $A_p(G)$}\}.$$
\begin{theorem}
Let $G$ be a discrete group. Then
\begin{enumerate}[label=({\roman*})]
\item $A_p(G)\subseteq  WCB_p(G)$.
\item Let $G$ be an amenable group. Then $WCB_p(G)=A_p(G)$.
\end{enumerate}
\end{theorem}
{\bf Proof.}
(i) If $x\in G$ and $\phi=\chi_x$, then $\Gamma_{\phi}(u)= u(x)\phi$. Hence $\Gamma_{\phi}$ is compact. Since the linear span of $\{\chi_x : x\in G\}$ is dense in $A_p(G)$. It follows that $A_p(G)\subseteq WCB_p(G)$. \\
(ii) Let $u\in WCB_p(G)$ and $(e_{\alpha})$ be a bounded approximate identity (B.A.I) for $A_p(G)$. Then there exist a subnet $(e_{\alpha\beta})$ and $v\in A_p(G)$ such that $e_{\alpha\beta}u\longrightarrow v$ in the weak topology of $A_p(G)$. Now, let $x\in G$ and $w\in A_p(G)$ such that $w(x) = 1$. Since $e_{\alpha\beta}(uw)\longrightarrow uw$ in the norm topology of $A_p(G)$. So,
$$v(x) = \lim_{\beta}(e_{\alpha\beta} u)(x) = \lim_{\beta}(e_{\alpha\beta} uw)(x) = w(x)u(x) = u(x).$$
Hence $u = v\in A_p(G)$. \\

{\bf Acknowledgment.} We would like to thank the Banach algebra center of Excellence for Mathematics Univercity of Isfahan.

\bibliographystyle{amsplain}

\begin{thebibliography}{5}
\bibitem{1}
M. Cowling, An application of Littlewood-Paley theory in harmonic analysis, Math. Ann., 241
(1979), 83.96.
\bibitem{2}
N. Dunford and J. T. Schwartz, Linear operators I, Interscience, New York (1958).
\bibitem{3}
P. Eymard, L'algebre de Fourier d.un groupe localement compact, Bull. Soc. Math. France, 92
(1964), 181.236.
\bibitem{4}
B. Forrest, Arens regularity and the $A_p(G)$ algebras, Proc. Amer. Math. Soc., 119 (1993),
595-598.
\bibitem{5}
E. E. Granirer, Exposed points of convex sets and weak sequential convergence, Mem. Amer.
Math. Soc., 123 (1972).
\bibitem{6}
----, On some spaces of linear functionals on the algebras $A_p(G)$ for locally compact groups,
Colloq. Math., LII (1987), 119132.
\bibitem{7}
C. Herz, Harmonic synthesis for subgroups, Ann. Inst. Fourier (Grenoble), 23 (1973), 91.123.
\bibitem{8}
A. T. Lau and A. \"Ulger, Some geometric properties on the Fourier Stieltjes algebras of locally
compact groups, Arens regularity and related problems, Trans. Amer. Math. Soc., 337 (1993),
321.359.
\end{thebibliography}

\authorsaddresses{ Department of Mathematics, University of Isfahan, 81745 Isfahan,Iran 
E-mail address: ghaeidamini@mau.ac.ir} 

\authorsaddresses{Department of Mathematics, University of Isfahan, 81745 Isfahan,Iran 
E-mail address: rejali@sci.ui.ac.ir}

\end{document}